\documentclass{article}
\usepackage{amsmath}
\usepackage{amscd}
\usepackage{latexsym}
\usepackage{amssymb}
\usepackage{amsmath}
\usepackage{hyperref}
\renewcommand{\abstractname}

\input xy \xyoption{all}

\date{}
\title{On some constructions of nil-clean, clean, and exchange rings}
\author{Alin Stancu}

\begin{document}
\maketitle
\newtheorem{theorem}{Theorem}
\newtheorem{corollary}{Corollary}
\newtheorem{lemma}{Lemma}
\newtheorem{proposition}{Proposition}
\newtheorem{example}{Example}
\newtheorem{rem}{Remark}
\newcommand{\G}{\Gamma}
\newcommand{\g}{\gamma}
\vspace{-7mm}
\date{}
{\noindent ${}^1$\textit{Department of Mathematics, Columbus State University, Columbus, GA 31907, \linebreak[0]
U.S.A.,\linebreak[0]
stancu\_alin1@columbusstate.edu}\\
}

\abstractname{}
\begin{abstract}\noindent
In this paper we discuss several constructions that lead to new examples of nil-clean, clean, and exchange rings.
A characterization of the idempotents in the algebra defined by a 2-cocycle is given and used to prove some of the algebra's properties (the
infinitesimal deformation case). From infinitesimal deformations we go to full deformations
and prove that any formal deformation of a clean (exchange) ring is itself clean (exchange).  Examples of nil-clean, clean, and exchange rings
arising from poset algebras are also discussed.

\end{abstract}

\section{Preliminaries}
Throughout the paper $k$ will be a commutative ring with unit and $A$ will be an associative $k$-algebra with unit. Bimodules $M$ over $A$ will be assumed
to be symmetric over $k$. That is, $am=ma$ for $a\in k$ and $m\in M$. We say that $A$ is $\mathbf{nil}$-$\mathbf{clean}$ if for each $a\in A$ we have a decomposition $a=e+n$ where $e$ is idempotent and $n$ is nilpotent in $A$. If for each $a\in A$ this decomposition is unique then we say that $A$ is $\mathbf{uniquely}$ $\mathbf{nil}$-$\mathbf{clean}$. We say that $A$ is $\mathbf{clean}$ if for each
$a\in A$ we have a decomposition $a=e+u$ where $e$ is idempotent and $u$ is invertible in $A$. If for each $a\in A$ this decomposition is unique then we say that $A$ is $\mathbf{uniquely}$ $\mathbf{clean}$. We say that $A$ is an $\mathbf{exchange}$ $\mathbf{ring}$ if for every $a\in A$ there exist an idempotent $e\in A$ such that $e\in aA$ and $1-e\in (1-a)A$. This is not the original definition of a ring with the exchange property but it is an equivalent characterization discovered independently by W. K. Nicholson and K. R. Goodearl and for convenience we will use it as our definition. We should note that this  condition is equivalent with the statement that for each $a\in A$ there exist an idempotent $e\in A$ such that $e\in Aa$ and $1-e\in A(1-a)$ as was shown by independently by R. B. Warfield and W. K. Nicholson (cf.\cite{TU}).

\section{Some background}  For an $A$-bimodule $M$ and $n\geq 0$ the $k$-module of the $k$-relative Hochschild $n$-cochains of $A$ with coefficients in $M$, denoted $C^n(A, M)$, is defined as follows: $C^0(A, M)=M$. If $n>0$ then a cochain $f\in C^n(A, M)$ is a $k$-multilinear function of $n$ variables $f:A\times\dots\times A\rightarrow M$.
The Hochschild coboundary map is the $k$-module map   $\delta^n:C^n(A, M)\rightarrow C^{n+1}(A, M)$ defined by  $(\delta^nf)(a_1,\dots,a_{n+1})=$ $a_1f(a_2,\dots,a_{n+1})+\sum(-1)^if(\dots,a_ia_{i+1},\dots)+(-1)^{n+1}f(a_1,\dots,a_n)a_{n+1}.$ It is easy to check that $\delta^{n+1}\delta^n=0$, so these maps form a cochain complex $C^\bullet(A, M)$. The kernel of $\delta^n$ is the $k$-module of the $k$-relative $n$-cocycles $Z^n(A, M)$ while the image is the $k$-module of relative $(n+1)$-coboundaries $B^{n+1}(A, M)$. Setting $B^0(A, M)=0$ the $k$-relative Hochschild cohomology modules of $A$ with coefficients in $M$ are defined by $H^n(A, M)=Z^n(A, M)/B^n(A, M)$. The word ``relative" is a reminder that Hochschild cohomology is a relative theory since $k$ is not necessarily a field. Because 2-cocycles will be used in this paper we should note that a $2$-cochain $f$ is a $2$-cocycle if for all $a, b, c\in A$ we have  \begin{equation}\label{cocycle} af(b,c)-f(ab,c)+f(a,bc)-f(a,b)c=0.\end{equation}
Note that if $a, b, c$ are all equal to some idempotent $e\in A$ then we have
\begin{equation}\label{idemp}
 ef(e, e)-f(e, e)e=0.
\end{equation}
Note also that if $d\in A$  and we take $a=d$ and $b=c=1$ then we get \begin{equation}\label{inv} f(d, 1)=df(1, 1). \end{equation}
Similarly, for $a=b=1$ and $c=d$  we  have \begin{equation}\label{inv2} f(1, d)=f(1, 1)d. \end{equation}

It is a classical result that $H^2(A, M)$ classifies the equivalence classes of $k$-allowable singular extensions of $A$ by $M$ (cf. \cite{GS}). A singular extension of $A$ by $M$ is a short exact sequence $$\mathcal{E}:\xymatrix{0\ar[r]&M\ar[r]& B\ar[r]&A\ar[r] &0}$$
where $B$ is a $k$-algebra equipped with the structure map $k\longrightarrow B$, the $k$-algebra structure on $A$ is the same as the one defined by the composition $k\longrightarrow B\longrightarrow A$, and the $B$-bimodule structure acquired by $M$  as the kernel of $B\longrightarrow A$ is identical with the one obtained through the map $B\longrightarrow A$. Two singular extensions of $A$ by $M$ are called equivalent if there is a commutative diagram $$\xymatrix{\mathcal{E}:0\ar[r]&M\ar[r]\ar[d]^{id_M} & B\ar[r]\ar[d] &A\ar[r]\ar[d]^{id_A} &0\\
\mathcal{E'}:0\ar[r] & M\ar[r] &B'\ar[r]& A\ar[r]&0}$$
in which the map $B\longrightarrow B'$  is an algebra map whose composite with $k\longrightarrow B$ is equal to $k\longrightarrow B'$. The set of equivalence classes of singular extensions has a natural $k$-module structure (see for example \cite{GS}).
A singular extension is called $k$-allowable if the map $B\longrightarrow A$ splits as a map of $k$-bimodules. Any extension equivalent to a $k$-allowable extension is allowable and the equivalence classes of such extensions forms a submodule of the module of all extensions (cf. \cite{GS}).

If $\mathcal{E}$ is such an extension then $A$ can be regarded as a submodule of $B$ via the splitting map $A\longrightarrow B$ and $B$ is the direct sum of
the $k$-bimodules $A$ and $M$. With $B=A+M$ the multiplication on $B$ has the form \begin{equation} (a, m)(a', m')=(aa', am'+ma'+f(a, a'))\end{equation}
where $f\in C^2(A, M)$. The associativity of $B$ implies that $f$ is a cocycle and, conversely, if $f$ is a 2-cocycle then the multiplication defined above
makes $A+M$ into an algebra and defines an allowable extension.
\section{Idempotents and extensions}
\begin{theorem}\label{Prop1}
Let $A$ be an associative $k$-algebra, $M$ an $A$-bimodule, and $f$ be a 2-cocycle with coefficients in $M$. Let $B$ denote the $k$-algebra
corresponding to the $k$-allowable singular extension defined by $f$. Then:

1) An element $(e, t)$ is an idempotent in $B$ if and only if $e$ is an idempotent in $A$ and \begin{equation}\label{idemp2} et+te+f(e, e)=t. \end{equation}

2) If $e$ is an idempotent of $A$ and $x\in M$ then equation \ref{idemp2} has the solution \begin{equation} t=(1-2e)f(e, e)+ex-xe. \end{equation}

3) Equation \ref{idemp2} has a unique solution if and only if $ex=xe$ for all $x\in M.$ \end{theorem}

\noindent\textsc{Proof.} 1) The element $(e, t)\in B$ is an idempotent if and only if we have $$(e^2, et+te+f(e, e))=(e, t).$$ This means that $e$ is an idempotent in $A$ and $t$ satisfies (\ref{idemp2}).

2) Given an idempotent $e$, equation (\ref{idemp2}) is easy to solve, for $t$, if  $e$ commutes with $M$. That is to say, if $ex=xe$ for all $x\in M$. Indeed, in this case, equation (\ref{idemp2}) is equivalent to \begin{equation} f(e, e)=(1-2e)t. \end{equation} Since $(1-2e)^2=1$ we get that $t=(1-2e)f(e, e)$ is the unique solution of $(\ref{idemp2})$.

What is interesting to note is that $t=(1-2e)f(e, e)$ is a solution of (\ref{idemp2}) even if $e$ does not commute with $M$. Indeed, because of (\ref{idemp}) we have \begin{center} $et+te+f(e, e)=e(1-2e)f(e, e)+(1-2e)f(e, e)e+f(e, e)=$ \end{center} \begin{center} $ =-ef(e, e)-ef(e, e)+f(e, e)=-2ef(e, e)+f(e, e)=t.$ \end{center}
In addition, for $x\in M$ the element $t=ex-xe$ satisfies the equation \begin{equation}\label{zero} et+te=t, \end{equation}
so the elements of the form  $t=(1-2e)f(e, e)+ex-xe$ are solutions of (\ref{idemp2}).

3) It is clear that if \ref{idemp2} has a unique solution then $ex=xe$, for al $x\in M.$ We also showed earlier that
if $ex=xe$ for all $x\in M$ then the equation has the unique solution $t=(1-2e)f(e, e).$ Therefore, a unique solution occurs if and
only if $ex=xe$ for all $x\in M.$ $\Box$\medskip

A consequence of theorem \ref{Prop1} is that it tells us explicitly how the idempotents of $A$ can be lifted
to idempotents of $B$. Taking this into account we can prove proposition \ref{Theorem 1}.  Note that
H. Lin and C. Xia gave in \cite{LX} a different proof for part 2) of proposition \ref{Theorem 1} and A. Diesl, in \cite{Diesl}, proved part 1) of \ref{Theorem 1} for the trivial extension.

\begin{proposition}\label{Theorem 1}

Let $A$ be an associative $k$-algebra, $M$ an $A$-bimodule, and $f$ a 2-cocycle with coefficients in $M$. Let $B$ denote the $k$-algebra
corresponding to the $k$-allowable singular extension defined by $f$. Then:

1) $B$ is nil-clean if and only if $A$ is nil-clean.

2) $B$ is clean if and only if $A$ is clean.

3) $B$ is exchange if and only if $A$ is exchange.

 \end{proposition}

  \noindent\textsc{Proof.}

  Clearly if $B$ is nil-clean (clean, exchange) then $A$ is nil-clean (clean, exchange) by the virtue of being a homomorphic image of $B$, so the interesting
  problem is going from $A$ to $B$.

To prove 1) lets assume that $A$ is nil-clean and let $(a, m)\in B$. Then there is an idempotent $e\in A$  and a nilpotent $x\in A$
such that $a=e+x$. This implies that  we have \begin{equation}\label{ne} (a, m)=(e, (1-2e)f(e, e))+(x, m-(1-2e)f(e, e)). \end{equation}
Since $(e, (1-2e)f(e, e)))$ is an idempotent it remains to show  that $(x, m-(1-2e)f(e, e))$ is nilpotent. It is an easy exercise to see that
 $(x, p)$ is a nilpotent in $B$ if and only if $x$ is  nilpotent in $A$. This is because if $x^n=0$ for some positive integer $n$ then we have $(x, p)^{2n}=(0, 0)$.

 To prove 2) lets assume that $A$ is clean and let $(a, m)\in B$. Then there is an idempotent $e\in A$  and an invertible element
 $u\in A$ such that $a=e+u$. Consider the decomposition  \begin{equation}\label{ce} (a, m)=(e, (1-2e)f(e, e))+(u, m-(1-2e)f(e, e)). \end{equation}

We will now prove that if $d$ is invertible in $A$ then $(d, p)$ is invertible in $B$ for all $p\in M$. (The converse is also true, so this characterizes the invertible elements in $B$.) To see this we construct the inverse of $(d, p)$. In particular this shows that $(u, m-(1-2e)f(e, e))$ is invertible in $B$, so it follows that $B$ is clean.

The element $(d^{-1}, -d^{-1}pd^{-1}-d^{-1}f(d, d^{-1})-d^{-1}f(1, 1))$ is the inverse of $(d, p)$. Indeed, we have \begin{center} $(d, p)(d^{-1}, -d^{-1}pd^{-1}-d^{-1}f(d, d^{-1})-d^{-1}f(1, 1))=$\end{center} \begin{equation}=(1, -pd^{-1}-f(d, d^{-1})-f(1, 1)+pd^{-1}+f(d, d^{-1})=(1, -f(1, 1)).\end{equation}
We also have \begin{center} $(d^{-1}, -d^{-1}pd^{-1}-d^{-1}f(d, d^{-1})-d^{-1}f(1, 1))(d, p)= $\end{center} \begin{center} $=(1, d^{-1}p-d^{-1}p-d^{-1}f(d, d^{-1})d-d^{-1}f(1, 1)d+f(d^{-1}, d)=$\end{center} \begin{equation}=(1, -d^{-1}f(d, d^{-1})d-d^{-1}f(1, 1)d+f(d^{-1}, d)). \end{equation}
We need to prove that \begin{equation}\label{inv3} -d^{-1}f(d, d^{-1})d-d^{-1}f(1, 1)d+f(d^{-1}, d)=-f(1, 1). \end{equation}
Recall that $f$ is a cocycle so for any $a, b, c\in A$ we have $$af(b,c)-f(ab,c)+f(a,bc)-f(a,b)c=0.$$ If $a=c=d^{-1}$ and $b=d$ then we get
\begin{equation}\label{13} d^{-1}f(d, d^{-1})-f(1, d^{-1})+f(d^{-1}, 1)-f(d^{-1}, d)d^{-1}=0.\end{equation}  Multiplying (\ref{13}) by $-d$ on the right we obtain
\begin{equation} \label{14} -d^{-1}f(d, d^{-1})d+f(1, d^{-1})d-f(d^{-1}, 1)d+f(d^{-1}, d)=0.\end{equation} Using (\ref{inv}) and (\ref{inv2}) in equation (\ref{14}) we obtain $$-d^{1}f(d, d^{-1})d+f(1, 1)-d^{-1}f(1, 1)d+f(d^{-1}, d)=0,$$ which proves (\ref{inv3}).

3) Assume now that $A$ is exchange. We use the following characterization of exchange rings. Namely, a ring $R$ is an exchange ring if and only if there exist an
ideal $I\subset J(R)$ such that $R/I$ is an exchange ring and all idempotents of $R/I$ can be lifted to idempotents of $R$. (cf. \cite{TU}, theorem 2.10)
We proved that each idempotent $e$ of $A$ lifts to at least one idempotent of $B$, namely $(e, (1-2e)f(e, e))$. We also have that $M^2=0$ and  $B/M$ is ring isomorphic to $A$. These two results combined imply that $B$ is exchange. $\Box$\medskip

$\mathbf{Remark.}$ One may try to show that $B$ is exchange by using the characterization stated in section 1 as the definition of exchange rings. That is, given $(a, m)\in B$, one has to show that there is an idempotent $(e, t)$ and elements $(b, c)$ and $(p, q)$ in $B$ such that $(e, t)=(a, m)(b, c)$ and $(1-e, -f(1, 1)-t)=(1-a, -f(1, 1)-m)(p, q)$. Since $A$ is exchange, for every $a\in A$ there exist an idempotent $e\in A$ and $r, s\in A$ such that $e=ar$ and $1-e=(1-a)s$. Recall that for each $x\in M$ the elements of the form $(e, (1-2e)f(e, e)+ex-xe)$ are idempotents of $B$. We have that
\begin{center} $(e, (1-2e)f(e, e))-(a, m)(re, f(r, e)-2rf(e, e)-rf(a, r)-rmr)=$\end{center} \begin{equation} \label{exchange} =(0, f(e, e)-af(r, e)+ef(a, r)+emr-mre-f(a, re)). \end{equation}
Since $f$ is a cocycle we have \begin{center} $af(r, e)-f(e, e)+f(a, re)-f(a, r)e=0$,\end{center} so (\ref{exchange}) is equal to $(0, ef(a, r)-f(a, r)e+emr-mre)$. Thus, for $x=-f(a, r)-mr$, we obtain \begin{center} $(e, (1-2e)f(e, e)+ex-xe)=(a, m)(re, f(r, e)-2rf(e, e)-rf(a, r)-rmr).$ \end{center}
Presumably one would now try to see that, for $t=(1-2e)f(e, e)+ex-xe$, there are $p\in A$ and $q\in M$ such that $(1-e, -f(1, 1)-t)=(1-a, -f(1, 1)-m)(p, q)$.
We don't know if this is actually true. Note that $t$ defined above depends on $r$, so if the answer is yes there should be some relation between $r$ and $s$ via the cocycle $f$. Whatever the answer is, the computations do not seem very pleasant. However, this raises the following more general question: Is it possible to find an exchange ring $R$ and elements $a, e\in R$, $e^2=e$ and $e\in aR$, but $1-e\notin(1-a)R$?

As a consequence of proposition \ref{Theorem 1} we obtain the next corollary. Using a different approach W. K. Nicholson proved in \cite{NZ}, for $ideal$-$extensions$, a similar result as part 2) of our corollary.

\begin{corollary}
Let $A$ be an associative $k$-algebra, $M$ an $A$-bimodule, and $f$ a 2-cocycle. Let $B$ denote the $k$-algebra
corresponding to the $k$-allowable singular extension defined by $f$. Then:

1) $B$ is uniquely nil-clean if and only if $A$ is uniquely nil-clean and the idempotents of $A$ commute
with all elements of $M$.

2) $B$ is uniquely clean if and only if $A$ is uniquely clean and the idempotents of $A$ commute
with all elements of $M$.

\end{corollary}

\noindent\textsc{Proof.}
Suppose that $B$ is uniquely nil-clean. If we assume that $A$ is not uniquely nil-clean then there is an $a\in A$ such that $a=e+x=e'+x'$ where $e$ and $e'$ are distinct idempotents and $x$ and $x'$ nilpotents. This  produces two distinct decompositions, corresponding to $e$ and $e'$, as in (\ref{ne}) and contradicts the assumption that $B$ is uniquely nil-clean. Also, since $B$ is uniquely nil-clean, for each idempotent $e\in A$ equation (\ref{idemp2}) has a unique solution, so $e$ commutes with all elements of $M$.

Conversely, assume that $A$ is uniquely nil-clean and $em=me$ for all $m\in M$ and all idempotents $e\in A$. Suppose that we have decompositions  \begin{center} $(a, m)=(e', t')+(x', m-t')=(e'', t'')+(x'', m-t'')$ \end{center}
for some $(a, m)\in B$, where $(e', t')$, $(e'', t'')$ are idempotents of $B$ and $(x', m-t')$, $(x'', m-t'')$ are nilpotents of $B$. Because $e'$ and $e''$ are idempotents and $x'$ and $x''$ are nilpotents and $A$ is uniquely nil-clean we get that  $e'=e''$ and $x'=x''$. If we denote by $e$ the element $e'$ (and $e''$) then $t'$ and $t''$ are solutions of equation (\ref{idemp2}). On the other hand since the idempotents
of $A$ commute with the elements of $M$, equation (\ref{idemp2}) has the unique solution $t=(1-2e)f(e, e)$, so $t'=t''=t$ and $B$ is uniquely nil-clean.
An identical argument can be employed to prove 2). $\Box$\medskip

\section{Deformations of  rings}\label{deformations}
In this section we need some basic notions about the deformation theory of algebras. We adopt the terminology and notations used by M. Gerstenhaber and S. D. Schack in \cite{GS}. The reader should be aware that we will be concerned with deformations inside the category of associative algebras. It is a known fact that any deformation of a unital algebra is equivalent with one in which the unit is not changed (cf. \cite{GS}) and we will assume that this is the case for the deformations considered below.

Let $A$ be a $k$-algebra and let $\alpha$ denote its multiplication.
A formal deformation of the $k$-algebra $A$, in the category of associative algebras, is a $k[[t]]$-associative algebra
given by a $k[[t]]$-bilinear multiplication \begin{center} $\alpha_t:A[[t]]\times A[[t]]\rightarrow A[[t]]$ of the form $\alpha_t=\alpha+t\alpha_1+t^2\alpha_2+\dots,$ \end{center}
where each $\alpha_i$ is a $k$-bilinear map $A\times A\rightarrow A$ (extended to be $k[[t]]$-bilinear). We write $A_t$ for the deformed
algebra (i.e. the $k[[t]]$-module $A[[t]]$ with the multiplication $\alpha_t$.

The null deformation is the deformations for which $\alpha_i=0$ for all $i$. A deformation $A'_t$ is equivalent to $A_t$ if there exist a $k[[t]]$-algebra
isomorphism $f_t:A'_t\rightarrow A_t$ of the form $f_t=Id_A+tf_1+t^2f_2+\cdots,$ where each $f_i$ is a $k$-linear map $A\rightarrow A$ extended to be $k[[t]]$-linear. A deformation equivalent to the null deformation is called trivial. An algebra with only trivial deformations is called analytically rigid over $k$.

The definition above may be reformulated as follows: A deformation of the associative algebra $A$ is a $t$-adically complete, $t$-torsion
free $k[[t]]$-algebra $A_t$ equipped with a $k$-algebra isomorphism $A_t\otimes_{k[[t]]}k\rightarrow A$.

\begin{lemma} \label{lemma}
Let $A_t$ be a deformation of the associative algebra $A$ and let $\displaystyle f=\sum_{i=0}^\infty a_it^i$ be an element of $A_t$. Then $f$ is invertible in
$A_t$ if and only if $a_0$ is invertible in $A$.
\end{lemma}

\noindent\textsc{Proof.}
Clearly if $f$ is invertible in $A_t$ then $a_0$ is invertible in $A$ since $a_0=f(0)$. Assume now that $a_0$ is invertible in $A$. If $g=\displaystyle \sum_{i=0}^\infty b_it^i$ then we have \begin{equation}\label{idemp3} \alpha_t(f, g)=\sum_{k=0}^\infty\sum_{m+n+p=k}^\infty
\alpha_m(a_n, b_p)t^k,\end{equation}  where $\alpha_0=\alpha$.
We want to find, for each $i\geq 0$, coefficients $b_i$ such that the equations $a_0b_0=1$ and \begin{equation} \displaystyle\sum_{m+n+p=k}\alpha_m(a_n, b_p)=0  \end{equation} are satisfied for each $k\geq 1$.
This can be done inductively. If coefficients $b_0, b_1, \dots ,b_r$ are known then $b_{r+1}$ can be found from the equation corresponding to $k=r+1$,
\begin{equation} \sum_{m+n+p=r+1}\alpha_m(a_n, b_p)=0.\end{equation}
Note that this equation  is linear in $b_{r+1}$ because $b_{r+1}$ appears only in the term $\alpha_0(a_0, b_{r+1})=a_0b_{r+1}$ and since $a_0$ is invertible we can solve the equation for $b_{r+1}$.
Thus we exhibited a procedure that finds a right inverse for $f$. Similarly we can find a left inverse for $f$, say $h$.  Since $A_t$ is associative we have
\begin{center} $h=\alpha_t(h, 1)=\alpha_t(h, \alpha_t(f, g))=\alpha_t(\alpha_t(h, f), g)=\alpha_t(1, g)=g.$ $\Box$\medskip \end{center}

\begin{corollary}\label{radical} Let $A_t$ be a deformation of the associative $k$-algebra $A$. Then $t\in J(A_t)$.
\end{corollary}
\noindent\textsc{Proof.}
The Jacobson radical of $A_t$ is the intersection of all maximal right ideals, so it is enough to show that $t$ is
contained in each maximal right ideal. If $I$ is a maximal right ideal such that $t\notin I$ then the right ideal generated
by $I$ and $t$ is equal to $A_t$. Thus there is some $b\in A_t$ and $g\in I$ such that $g+tb=1$. This implies that $g=1-tb$, so $g$ is invertible
and we get that $I=A_t$. This is impossible, so we have $t\in I$ and consequently $t\in J(A_t)$.
$\Box$\medskip

The following lemma can be found in \cite{GS}. It is in some sense an extension of the lifting properties
for ideals that are nilpotent.

\begin{lemma} \label{lemma2}
Let $\hat{A}$ be an algebra and $I$ a two sided ideal of $\hat{A}$ such that $\hat{A}$ is complete
in the $I$-adic topology. Then any idempotent $e\in\hat{A}/I$ lifts to an idempotent $\hat{e}\in\hat{A}.$  In particular,
if $A_t$ is a deformation of $A$ and $e\in A$ is an idempotent then there is an idempotent $e_t\in A_t$ with constant term $e.$
\end{lemma}
\noindent\textsc{Proof.}
Following \cite{GS}, the idea is to construct a polynomial $f\in\mathbb{Z}[X]$ with the following property: If $a^2-a\in J$ for some ideal $J\subset\hat{A}$ then $f(a)-a\in J$ and $f(a)^2-f(a)\in J^2$. When $f$ is such a polynomial, starting with $a=e$ and $J=I$, we obtain a Cauchy sequence
$\{f^n(e)\}$ which converges to an idempotent $\hat{e}\in A_t$ that lifts $e$. The polynomial $f$ can be found by formally applying  Newton's method
to the equation $x^2-x=0$. Begin with the approximate solution $x_0=a$. We have $x_1=x_0-\frac{f(x_0)}{f'(x_0)}=\frac{-a^2}{2a-1}$ and $x_1^2-x_1=
\frac{(a^2-a))^2}{(2a-1)^2}$. Because an inverse for $(2a-1)$ modulo $J^2$ is $(2a-1)(4a^2-4a-1)$ the polynomial $f(X)=-X^2(2X-1)(4X^2-4X-1)$ has the required property. $\Box$\medskip

\begin{lemma}\label{lemma3} Let $A$ be a $k$-algebra and let $e\in A$ be central idempotent. If $A_t$ is a deformation of $A$ then
$e$ lifts uniquely to an idempotent $e_t\in A_t$.

\end{lemma}
\noindent\textsc{Proof.}
By lemma \ref{lemma2} we know that $e$ lifts to an idempotent $e_t\in A_t,$ say $e_t=e+a_1t+a_2t^2+\dots$. We will show that the coefficients $a_i$ are
uniquely determined by $e$ and the deformation cochains $\alpha_i$. To see this we use that $\alpha_t(e_t, e_t)=e_t,$ so for each $k\geq 1$ we have the equations \begin{equation}\label{idemp4} \displaystyle\sum_{m+n+p=k}\alpha_m(a_n, a_p)=a_k. \end{equation} For $k=1$ we get the equation $ea_1+a_1e+\alpha_1(e, e)=a_1.$ Because $e$ is central this equation has the unique solution $a_1=(1-2e)\alpha_1(e, e)$,
so $a_1$ is uniquely determined.
When $k=2$ we have that $a_2$ satisfies the equation $ea_2+a_2e+\beta_2(e, e)=a_2,$ where $\beta_2(e, e)=\alpha_0(a_1, a_1)+\alpha_1(a_1, e)+\alpha_1(e, a_1)+\alpha_2(e, e).$ This equation has the unique solution $a_2=(1-2e)\beta_2(e, e).$ We can now proceed by induction. Assume that $a_1, a_2,\dots, a_k$ are unique. Then we have that $a_{k+1}$ satisfies the equation $ea_{k+1}+a_{k+1}e+\beta_{k+1}(e, e)=a_{k+1},$ where $\displaystyle \beta_{k+1}(e, e)=\sum_{ m+n+p=k+1}\alpha_m(a_n, a_p)$ and $(n, p)\neq (0, k+1)$ and $(k+1, 0)$. This implies that $a_{k+1}=(1-2e)\beta_{k+1}(e, e),$ so $a_{k+1}$ is uniquely determined by $e$ and $\alpha_i$. $\Box$\medskip

$\mathbf{Remark}$ $\mathbf{1.}$ It is interesting to note the constraints if one tries to solve, for $a_i$, the equation $\alpha_t(f, f)=f,$ where $f=e+a_1t+a_2t^2+\dots$, without the assumption that $e$ is central. The first obstruction to ``integrating" $e$ to a full idempotent of $A_t$ is that $a_1$ should be a solution of the equation $ea_1+a_1e+\alpha_1(e, e)=a_1$,
which is nothing else than equation (\ref{idemp2}) for cocycle $\alpha_1$ (the associativity of $\alpha_t$ implies $\alpha_1$ is a cocycle in $Z^2(A, A)$.) We know that this equation has at least one solution, namely $a_1=(1-2e)\alpha_1(e, e)$, so $e$ is ``infinitesimally" integrable. For this choice of $a_1$ the second obstruction is that $a_2$ should be a solution of equation $ea_2+a_2e+\beta_2(e, e)=a_2,$ where $\beta_2(e, e)=\alpha_0(a_1, a_1)+\alpha_1(a_1, e)+\alpha_1(e, a_1)+\alpha_2(e, e).$ It can be shown that $e\beta_2(e, e)=\beta_2(e, e)e$, so $a_2=(1-2e)\beta_2(e, e)$ satisfies the second obstruction.
However, we don't know if  $a_i=(1-2e)\beta_i(e, e)$, for $i\geq 3$, are solutions of the equations $ea_{i}+a_{i}e+\beta_{i}(e, e)=a_{i},$ where $\displaystyle \beta_{i}(e, e)=\sum_{m+n+p=i}\alpha_m(a_n, a_p)$ with $(n, p)\neq (0, k+1)$ and $(k+1, 0)$. This would be true if $e\beta_i(e, e)=\beta_i(e, e)e$, but it is not clear that this is the case.

$\mathbf{2.}$ If $e$ is  not a central idempotent in $A$ then it is possible to have two lifts in $A_t$. To see this consider the case of the trivial deformation and let $x\in A$  such that $ex-xe\neq 0.$ Then $e$ can be viewed as an idempotent in $A[[t]]$ and one can also check that it has a nontrivial lifting to an idempotent $e_t=e+a_1t+(1-2e)a_1^2t^2+2(1-2e)(a_1^3-a_1^3e-ea_1^3)t^3+\dots,$ where $a_1=ex-xe$.

To prove our next theorem we need the following result due to J. Han and W. K. Nicholson (cf.\cite{N}). The proof
is short and for the clarity of our paper we include it.

\begin{proposition} \label{clean}
Let $I$ be an ideal of $A$ such that $I\subset J(A)$. Then $A$ is clean if and only if
$A/I$ is clean and idempotents lift modulo $I$.
\end{proposition}

\noindent\textsc{Proof.}
If $A$ is clean then $A/I$ is clean because it is a factor of $A$. Let $\overline{a}$ be an idempotent in $A/I$.
This implies that we have $a^2-a\in I$. Since $A$ is clean there is an idempotent $e\in A$ and an unit $u\in A$
such that $a=e+u$. Then $a-u^{-1}(1-e)u=u^{-1}(a^2-a)\in I,$ so $a$ lifts to the idempotent $u^{-1}(1-e)u$.
To prove the converse let $a\in A$. Since $A/I$ is clean there is an idempotent $\overline{e}\in A/I$ and an unit $\overline{u}\in A/I$
such that $\overline{a}=\overline{e}+\overline{u}$. We may assume that $e$ is an idempotent in $A$ because all idempotents of $A/I$ lift modulo $I$.
Thus we have that $a-e$ is a unit in $A/I$. This implies that $a-e$ is a unit in $A$ since $I\subset J(A)$.

\begin{theorem}\label{ce} Let $A$ be a $k$-algebra and $A_t$ a deformation of $A$. Then

1) $A$ is clean if and only if $A_t$ is clean.

2) $A$ is exchange if and only if $A_t$ is exchange.

\end{theorem}

\noindent\textsc{Proof.}
1) Clearly if $A_t$ is clean then $A$ is clean because it is a factor of $A_t$. To see that the converse is true note that $A/(t)$
is isomorphic to $A$, so it is clean. In addition, we have $(t)\subset J(A_t)$ (cf. corollary \ref{radical}) and all idempotents of $A$ lift to $A_t$ (cf. lemma \ref{lemma2}).
Therefore $A_t$ is clean (cf. proposition \ref{clean})

2) It follows easily as a consequence of theorem 2.10 in \cite{TU} which asserts that a ring $R$ is exchange if and only if there exist an
ideal $I\subset J(R)$ such that $R/I$ is an exchange ring and all idempotents of $R/I$ can be lifted to idempotents of $R$.
Combining this result and lemma \ref{lemma2} the proof is clear. $\Box$\medskip

\begin{corollary}\label{uniquelyclean} Let $A$ be a uniquely clean $k$-algebra. Then any deformation $A_t$ of $A$ is uniquely clean. \end{corollary}

\noindent\textsc{Proof.} Because $A$ is clean we get that $A_t$ is clean (cf. theorem \ref{ce}). Suppose that we have two decompositions in $A_t$, say  $h_t=e_t+u_t=e'_t+u_t'$ with $e_t, e'_t$ idempotents and $u_t, u_t'$ units in $A_t$. This implies that we have two decompositions of $h_t(0)$ in $A$, in sums of idempotents and units, namely $h_t(0)=e_t(0)+u_t(0)=e'_t(0)+u'_t(0)$. Since $A$ is uniquely clean we get $e_t(0)=e'_t(0)$. Because in a uniquely clean ring all idempotents are central we have that $e_t=e'_t.$ (cf. lemma \ref{lemma3}) . Therefore $A_t$ is uniquely clean. $\Box$\medskip

 $\mathbf{Remark}$ $\mathbf{1.}$ If $A_t$ is the trivial deformation of $A$ then $A_t$ is isomorphic to the $k$-algebra of
power series in one variable $A[[t]]$. Therefore $A$ is clean (exchange) if and only if $A[[t]]$ is clean (exchange) (cf. theorem \ref{ce}).

$\mathbf{2.}$ If $A$ is uniquely clean then $A[[t]]$ is uniquely clean (cf. corollary \ref{uniquelyclean}).
Clearly the converse is also true.

\section{Posets of Algebras}
Let $\mathcal{C}$ be a finite poset, say with $n$ elements,  viewed as a category in the usual way: for each $i\leq j$ there is a unique map $\varphi^{ij}:i\rightarrow j$.
Let $\mathbb{A}$ be a presheaf  of $k$-algebras over $\mathcal{C}$. That is: $\mathbb{A}$ is a functor $\mathbb{A}:\mathcal{C}^{op}\rightarrow k-\mathbf{alg}$. One may associate to each presheaf $\mathbb{A}$ a single algebra, denoted $\mathbb{A!}$.  The algebra $\mathbb{A}$ consists of $\mathcal{C}\times\mathcal{C}$ matrices $(a_{ij})$ with $a_{ij}\in\mathbb{A}(i)$ if $i\leq j$ and $a_{ij}=0$ otherwise. For simplicity we will denote $\mathbb{A}(i)$ by $\mathbb{A}^i$. The addition is componentwise and the multiplication $(a_{ij})(b_{ij})=(c_{ij})$ is induced by the matrix multiplication with the understanding that for $h\leq i\leq j$ the summand $a_{hi}b_{ij}$ of $c_{hj}$ is equal to $a_{hi}\varphi^{hi}(b_{ij})$. (The construction can be adapted to infinite posets and is due to M. Gerstenhaber and S. D. Schack (cf.\cite{GS}) and it plays an important role in studying the Yoneda cohomology of presheaves.) In this paper we will be concerned only with the case of a finite poset.
Lets look at an example. Consider a presheaf $\mathbb{A}$ of $k$-algebras over the poset given by the diagram below.
$$
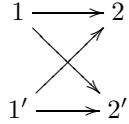
\xymatrix{1\ar[r]\ar[dr] & 2\\
1'\ar[ur]\ar[r] & 2'}$$
Corresponding to the the linear order $1, 1', 2, 2'$ (respectively $2, 1, 1', 2'$) we obtain that $\mathbb{A}!$ consists of matrices of the form

\begin{center}
$\left[                                                                                                                     \begin{array}{cccc}
                                                                                                                       $\;$ \mathbb{A}^1 & 0 &  $\;$ \mathbb{A}^1  &  $\,$ \mathbb{A}^1  \\
                                                                                                                       0 &  $\;$ $\;$ \mathbb{A}^{1'} &  $\;$ $\;$ \mathbb{A}^{1'} & $\;$ $\;$ \mathbb{A}^{1'} \\
                                                                                                                       0 & 0&  $\;$ \mathbb{A}^2  & 0 \\
                                                                                                                       0 & 0 & 0 &  $\;$ $\;$ \mathbb{A}^{2'}  \\
                                                                                                                     \end{array}
                                                                                                                   \right], \mathrm{respectively} \left[
                                                                                                                                           \begin{array}{cccc}
                                                                                                                                              $\;$ \mathbb{A}^2  &  0 & 0 &  0 \\
                                                                                                                                              $\;$ \mathbb{A}^1 &  $\;$ \mathbb{A}^1 &  0 &  $\;$ \mathbb{A}^1  \\
                                                                                                                                              $\;$ $\;$\mathbb{A}^{1'} &  0 & $\;$ $\;$ \mathbb{A}^{1'} & $\;$ $\;$ \mathbb{A}^{1'} \\
                                                                                                                                             0 & 0 & 0 &  $\;$ $\;$ \mathbb{A}^{2'}\\
                                                                                                                                           \end{array}
                                                                                                                                         \right].$
                                                                                                                   \end{center}

Note that different linear orders produce different representations of $\mathbb{A}!$ but, since the poset $\mathcal{C}$ is finite, one can always choose a representation of $\mathbb{A}!$ (in general not uniquely) in an upper triangular form. Taking this into account we can prove the following result:

\begin{theorem}\label{cleanposets}Let $\mathcal{C}$ be a finite poset and $\mathbb{A}$ a presheaf of $k$-algebras over $\mathcal{C}$. Then

1) $\mathbb{A}!$ is clean (respectively nil-clean) if and only if $\mathbb{A}^i$ is clean (respectively nil-clean) for all $i\in\mathcal{C}$.

2) $\mathbb{A}!$ is exchange if and only if $\mathbb{A}^i$ is exchange for all $i\in\mathcal{C}$.

\end{theorem}

\noindent\textsc{Proof.} 1) By taking a linear order that results in an upper triangular representation of $\mathbb{A}!$ we have that the strictly upper
triangular matrices form an ideal of $\mathbb{A}!$. If we denote this ideal by $I$ and assume that $\mathbb{A}!$ is clean (respectively nil-clean) then $\mathbb{A}!/I$
is clean (respectively nil-clean). This implies that the direct product $\prod_{i\in\mathcal{C}}\mathbb{A}_i$ is clean (respectively nil-clean), so each $\mathbb{A}_i$ is clean (respectively nil-clean). Conversely, if $\mathbb{A}_i$ are all clean (respectively nil-clean) and $M\in\mathbb{A}!$ then we can take decompositions of the elements on the main diagonal of $M$, $x_i=e_i+u_i$, $e_i$ idempotent and $u_i$ unit in $\mathbb{A}_i$ (respectively $x_i=e_i+n_i$, $e_i$ idempotent and $n_i$ nilpotent in $\mathbb{A}_i$) and write $M$ as the sum of the diagonal matrix $D$ with entries $e_i$ and the upper triangular matrix $M-D$. Clearly $D$ is idempotent in $\mathbb{A}!$, and $M-D$ is invertible in the clean case (respectively nilpotent in the nil-clean case).

2) We consider again a linear order that results in an upper triangular representation of $\mathbb{A}!$. If $\mathbb{A}!$ is exchange then any of its factors is also exchange. By taking $I$ to be the ideal of upper triangular matrices in $\mathbb{A}!$ we  get that $\mathbb{A}!/I\cong\prod_{i\in\mathcal{C}}\mathbb{A}^i,$ so each $\mathbb{A}^i$ is exchange. For the converse we should note that $I\subseteqq J(\mathbb{A}!)$ and each idempotent in $\mathbb{A}!/I$ lifts to an idempotent in $\mathbb{A}!$. $\Box$\medskip

$\mathbf{Example}$ $\mathbf{1.}$ Consider the poset given by the relations $2'\geq 1\leq 2$. Let $A=\mathbb{A}^1=\mathbb{Z}_2[X]/(X^2)$ and  $\mathbb{A}^2=\mathbb{A}^{2'}=\mathbb{Z}_2$. Then we can write $\mathbb{A}!$ as \begin{center} $\mathbb{A}!=\left[
                                                                                  \begin{array}{ccc}
                                                                                     A &  A& A \\
                                                                                    0 & $\;$ \mathbb{Z}_2 & 0 \\
                                                                                    0 & 0 & $\;$ \mathbb{Z}_2 \\
                                                                                  \end{array}
                                                                                \right]$\end{center}
Since both $A$ and $\mathbb{Z}_2$ are nil-clean rings we have that $\mathbb{A}!$ is a nil-clean ring. In particular, $\mathbb{A}!$ is clean. Actually it can be shown that $\mathbb{A}!$ is strongly clean.

$\mathbf{2.}$ Consider the presheaf given by the poset $\mathcal{C}$ below, where $\mathbb{A}^i=k$ for all $i\in \mathcal{C}$.
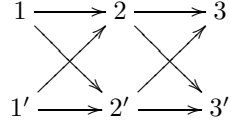
$$\xymatrix{1\ar[r]\ar[dr] & 2\ar[r]\ar[dr] &3\\
1'\ar[ur]\ar[r] & 2'\ar[ur]\ar[r] &3'}$$
If $k$ is clean (nil-clean, exchange) then $\mathbb{A}!$ is
clean (nil-clean, exchange) (cf. theorem\ref{cleanposets}).

The Hochschild cohomology of $\mathbb{A!}$ has the property that $H^2(\mathbb{A}!, \mathbb{A}!)=k$ and $H^3(\mathbb{A}!, \mathbb{A}!)=0$. To see this one should note that the 2-sphere is the geometric realization of the above poset. A rather special case of the $Cohomology$ $Comparison$ $Theorem$, of M. Gerstenhaber and S. D. Schack, then provides a natural isomorphism between the Hochschild cohomology of $\mathbb{A}!$ and the simplicial cohomology of the nerve of $\mathcal{C}$ (cf.\cite{GS}).

In particular, there are 2-cocycles in $C^2(\mathbb{A}!, \mathbb{A}!)$ which are not coboundaries ($H^2(\mathbb{A}!, \mathbb{A}!)\neq 0$). In addition,  every such cocycle can be integrated to a deformation since there are no obstructions ($H^3(\mathbb{A}!, \mathbb{A}!)=0).$  Finally, each such deformation is non-trivial (not isomorphic to the ring of power series over $\mathbb{A}!$) since the cocycle is not a coboundary. Therefore, if $k$ is clean (exchange) any such deformation of $\mathbb{A}!$ is not trivial and is clean (exchange).

\end{document}